\documentclass[12pt,leqno]{amsart}
\usepackage[dvips]{graphicx}
\usepackage{amsfonts}
\usepackage{amssymb}
\usepackage{amsmath}
\usepackage{amscd}
\usepackage{graphicx}
\def\DATE{\today}
\newtheorem{theorem}{Theorem}
\newtheorem{definition}[theorem]{Definition}
\newtheorem{corollary}[theorem]{Corollary}

\newtheorem{lemma}[theorem]{Lemma}

\newtheorem{proposition}[theorem]{Proposition}

\newcommand\C{\mathbb{C}}
\newcommand\R{\mathbb{R}}

\newcommand\ad{\mathrm{ad}}

\newcommand\g{\mathfrak{g}}
\newcommand\G{\Gamma}
\newcommand\h{\mathfrak{h}}

\newcommand\n{\mathfrak{n}}
\newcommand\K{\mathbb{K}}

\newcommand\A{\mathcal{A}}

\newcommand\pf{\noindent{\it Proof. }}
\newcommand\lr{\left\{ \begin{array}{l}}
\def\ds{\displaystyle}

\newcommand{\mg}{\mathfrak g }

\newcommand{\mk}{\mathfrak k }

\newcommand{\mgg}{\mathfrak g }

\newcommand{\ra}{\rightarrow}

\newcommand{\rmc}{\textrm{\rm C}}

\newcommand{\lra}{\longrightarrow}

\title{Symplectic structures on $2$-step nilpotent  Lie algebras }
\author{Elisabeth Remm, Michel Goze}
\date{28 Tichri 5776}
\address{E.R:Universit\'{e} de Haute Alsace, LMIA, 4 rue des Fr\`{e}res Lumi\`{e}%
re, F.68093 Mulhouse. \ \
M.G: Ramm Algebra Center, F68800 Rammersmatt}
\email{elisabeth.remm@uha.fr, goze.rac@gmail.com}

\begin{document}

\maketitle

\begin{abstract} We study symplectic structures on nilpotent Lie algebras. Since the classification of nilpotent Lie algebras in any dimension seems to be a crazy dream, we approach this study in case of $2$-step nilpotent Lie algebras (in this sub-case also, the classification fo the dimension greater than 8 seems very difficult), using not a classification but a description of subfamilies associated with the characteristic sequence. We begin with the dimension $8$, first step where the classification becomes difficult.

\end{abstract}

\noindent{\bf Introduction}

Let $M$ be a $2p$-dimensional  differentiable manifold. A symplectic form on $M$ is a differential $2$-form $\theta$ satisfying $\theta ^p=\theta \wedge \cdots \wedge \theta \neq 0$
and $d \theta =0$, where $d$ is the derivation operator of the De Rham cohomology of $M.$ The existence problem of a symplectic form on an even dimensional differentiable manifold  is very classical.
Although any point $x \in M$ has a neighborhood $V(x)$  admitting locally a symplectic form, it is often very difficult to extend it to the whole variety.  We get around this problem when $M$ is a Lie group $G$ and the symplectic form is left invariant. In this case the existence problem reduces to the existence on the Lie algebra $\frak{g}$ of the Lie group $G$ of an exterior $2$-form (antisymmetric bilinear form) of maximal rank $2p$ and closed for the cohomology of $\frak{g}$ with values in $\mathbb{R}.$ It will be similar when $M$ is a nilvariety $M=\Gamma  / G$ with $G$ is a nilpotent Lie group and $\Gamma$ a cocompact discret subgroup. In fact any symplectic form on $M$ is cohomologous to an invariant form on $G$. So the existence of a symplectic form on $M$ reduces to the existence of an exterior $2$-form on $\mathfrak{g}$. So we will focus on this last problem.

% Locally, all the symplectic forms have the same expression from the Darboux's theorem. Moreover, there exists always %in a local chart of $M$ a symplectic form, the existence problem consists to extend these local symplectic forms to the %full variety.   All cotangent bundles admit canonical symplectic forms, a fact relevant for analysis of differential %operators, dynamical systems, classical mechanics, etc.  In the particular case of a nilmanifold $M = \Gamma / G$, that %is, a compact quotient of a nilpotent simply connected Lie group by a cocompact discrete subgroup $\Gamma$,  any %symplectic form on $M$ is cohomologous to a left invariant form on $G$ and hence it is represented by a closed 2-form %on the Lie algebra $\g$ of $G$.
%Therefore, the study of existence of symplectic structures on the nilmanifold $\Gamma /G$ reduces to the existence of %a nondegenerate closed $2$-form on the Lie algebra $\g$, that is which is a $2$-cocycle for the cohomology of $\g$ %with values in $\R$.

To simplify, we will call symplectic Lie algebra an even dimensional  Lie algebra over the field $\K=\R$ or $\C$ which admits a symplectic (bilinear) form. The problem of determining whether an arbitrary
Lie algebra admits symplectic structures is difficult in general. It is solved for nilpotent Lie algebras up to the dimension $6$ (\cite{BG}) because the classification of nilpotent Lie algebras is known up to the dimension $7$. Several sub-families of nilpotent symplectic Lie algebras have been described. For example, symplectic filiform algebras, that is nilpotent Lie algebras with a maximal nilindex $n=\dim\g -1$, are studied in \cite{Bu, GJK, Mi}, symplectic Lie algebras which are of Heisenberg type in \cite{DT}, symplectic
nilpotent Lie algebras associated with graphs \cite{PO-TI},
or free nilpotent Lie algebras \cite{dB1}.

In this paper, we are interested by symplectic nilpotent Lie algebra which are $2$-step nilpotent that is with a nilindex equal to $2$. Recall  that the differential geometry of the $2$-step nilpotent Lie groups is the subject of many papers. Of course, we know sub-families of symplectic $2$-step nilpotent Lie algebras, as the nilradicals of parabolic subalgebras of split simple Lie algebras \cite{CadB}. We process here using the description of $2$-step nilpotent Lie algebras based on the characteristic sequence (see \cite{GRKegel}). As the dimension $8$ is the first that we meet without classification results, we are obliged to expose a method which is not based on the classification. Thus, in this paper, we are principally interested by this dimension.

\section{Symplectic structures on  Lie algebras}

\subsection{Generalities}

Let $\g$ be a $n$-dimensional $\K$-Lie algebra. We denote by $\g^*$ the dual vector space of $\g$ and by  $\Lambda^p\g^*$  the space of skew-symmetric $p$-linear forms on $\g$.  For any $\alpha \in \g^*$, we consider the skew-symmetric $2$-form $\theta=d\alpha$ which is defined as follows:
$$d\alpha(X,Y)=-\alpha[X,Y]$$
for any $X,Y \in \g$. For any $\theta \in \Lambda^2\g^*$, we define the  skew-symmetric $3$-form $d\theta$ given by
$$d\theta(X,Y,Z)=\theta([X,Y],Z)+\theta([Y,Z],X)+\theta([Z,X],Y)$$
for any $X,Y,Z \in \g$. In particular we have
$$d(d\alpha)=0$$
for any $\alpha \in \g^*$. The $2$-form $\theta$ is called closed if $d\theta =0$ and it is called exact if $\theta=d\alpha$ for some $\alpha \in \g^*.$  In fact, these operators $d$ coincide with the differential operators associated with the Chevalley-Eilenberg complex of $\g$:
\begin{equation}\label{eq:complex}
0\lra\K \lra \mgg^*\stackrel{d_1}{\lra}\Lambda^2\mgg^*\stackrel{d_2}{\lra}\ldots \ldots\stackrel{d_{n-1}}{\lra} \Lambda^n\mg^*\lra 0\;\,,
\end{equation}
with
\begin{equation} d_p \Omega\,(X_1,\ldots,X_{p+1})=\sum_{1\leq i<j\leq p+1}(-1)^{i+j-1}\Omega([X_i,X_j],X_1,\ldots,\hat{X_i},\ldots,\hat{X_j},\ldots,X_{p+1}). \label{eq:diff}\end{equation}
for any $\Omega \in \Lambda^p\g^*$ and
 $X_1,\ldots,X_{p+1} \in \g$.
The cohomology of this complex $(\Lambda^*\g^*,d)$  is denoted by $H^*(\g)$.

\begin{definition}
Let $\g$ be a $n=2p$-dimensional $\K$-Lie algebra.
A symplectic form on $\g$ is a closed 2-form $\theta$ which is also nondegenerate that is
$$d\theta =0, \ \ \theta^p \neq 0.$$
A Lie algebra provided with a symplectic form $\theta$ will be called a symplectic Lie algebra and denoted by the pair $(\g,\theta)$.
\end{definition}

A particular case corresponds to the Frobeniusian Lie algebras. Assume that $\g$ is a symplectic Lie algebra with a symplectic form $\theta$ which is exact, that is there is $\alpha \in \g^*$ such that $\theta = d\alpha$.
In this case the symplectic Lie algebra, $(\g,\theta)=(\g,d\alpha)$ is called Frobeniusian. Let us note that the linear form $\alpha$ is, in this case, of maximal Cartan class equal to $2p=\dim \g$. We deduce immediately the well known property that
we have an open orbit in the coadjoint representation.

\medskip

Considering a symplectic Lie algebra, we can be interested by additional geometrical structures adapted to the symplectic form such as a pseudo-riemannian structure or a complex structure (or both simultaneously) and also by geometrical properties associated with the symplectic one such as an affine structure.

\begin{definition}
Let $\g$ be a $n$-dimensional $\K$-Lie algebra. An affine structure on $\g$  is a $\K$-bilinear product $\g\times \g \rightarrow \g$ denoted $XY$
satisfying
\begin{equation}
\label{ls}
\left\{
\begin{array}{l}
   XY-YX=[X,Y],\\
   X(YZ)-(XY)Z=Y(XZ)-(YX)Z,
\end{array}
\right.
\end{equation}
 for any $X,Y,Z \in \g.$
\end{definition}
A Lie algebra over $\K$ admitting an affine structure is also called affine. Let us note that the nonassociative product $XY$ is Lie-admissible and often called left-symmetric. In general it is very difficult for a given Lie algebra to decide whether it is affine or not. But we know partial results. For example, it is known that a $\K$-Lie algebra $\frak{g}$  satisfying $\g = [\g, \g]$ is not affine. Likewise,  any complex nilpotent Lie algebras of dimension $n \leq7$ is affine. But there exist  examples of nilpotent Lie algebras of dimension $10$ which are not affine.

\begin{proposition}
Any symplectic Lie algebra is affine.
\end{proposition}
\pf Let $(\g,\theta)$ be a symplectic Lie algebra. We consider the product $XY$ given by
$$XY=f(X)Y$$
where $f: \g \rightarrow End(\g)$ is defined implicitly by
$$\theta (f(X)(Y),Z)=-\theta (Y,[X,Z])$$
for any $X,Y,Z \in \g$. Since $\theta$ is symplectic, this product $XY$ is well defined and provides $\g$ with an affine structure. 

Let $\{ X_i  \}$ be a basis of $\frak{g}$ for which the symplectic form $\theta$ is reduced to the canonical form. Let $A$ be its matrix
$A=I_p\otimes S$ with
$S=\left(
\begin{array}{rr}
0 & 1 \\
-1 & 0
\end{array}
\right) .$
We denote by $T_X$ the matrix of $ad_X$ with respect to the basis $\{ X_i \}$ and by $M_X$ the matrix of $f(X).$ We have
$$^t (M_X Y)AZ=-^tYAT_X Z,$$
this gives $^tM_XA=-AT_X.$
But $A^{-1}=-A=^t\! \! A.$ So $^tM_X=AT_XA$ and
$$M_X= A \, ^tT_X A.$$

When $\theta$ is a frobeniusian form, that is $\theta=d\alpha$, the associated affine structure is written
$$\alpha [XY,Z]=-\alpha[Y,[X,Z]].$$

\medskip

As we have said above, we can consider additional geometrical structures on a symplectic Lie algebra. One of the most interesting is the  quadratic structure.
\begin{definition}
Let $\g$ be an $n$-dimensional $\K$-Lie algebra.
We say that $\g$ is a quadratic Lie algebra if there exists  a non-degenerate symmetric bilinear form $B$ on $\g$ such that $$B([x, y], z) = B(x, [y, z])$$
 for all $x, y, z \in \g$. In this case, we will say that $B$ is an invariant scalar product on $\g$ and we will denote by $(\g,B)$ this quadratic Lie algebra.
 \end{definition}

 If $\g$ is also a symplectic Lie algebra (and $n=2p$), we have to consider the triple $(\g,\theta,B)$ and we will speak about symplectic quadratic Lie algebra. 
 \begin{proposition}
If $(\g,\theta,B)$ is a symplectic quadratic Lie algebra then there exists a $B$-skew-symmetric invertible derivation $D$  of $\g$
such that
$$\theta(X,Y)=B(D(X),Y)$$
for all $X,Y \in \g$.
\end{proposition}
Conversely, if $(\g,B)$ is a quadratic Lie algebra and $D$ an inversible derivation, then $\theta$ defined by the previous formula is symplectic as soon as $n=2p$.
The characterization of nilpotent quadratic symplectic Lie algebras is due to Benson and Gordon (\cite{BE-GO}).

To end these remarks, we can consider also complex structures on symplectic Lie algebras.

\begin{definition}A  complex structure on an $n=2p$-dimensional $\R$-Lie algebra $\g$ is a linear endomorphism $J$ of $\g$ such that:
\begin{enumerate}
  \item $J^2 = - Id$;
  \item $ [JX,JY] - [X,Y] - J[JX,Y] - J[X,JY] =0, \qquad \forall X, Y \in \g$.
\end{enumerate}
\end{definition}
Assume now that  $(\g,\theta)$ is a symplectic Lie algebra. A K\"ahler structure  on $(\g,\theta)$ is given by a complex structure $J$, which  is compatible with $\theta$, that  is
$$\theta (JX, JY) = \theta (X,Y) \qquad \text{ for all } X, Y \in \g.$$
Any such K\"ahler structure defines a pseudo-Riemannian metric $B$ on the Lie algebra in the following form:
$$B(X,Y) = \theta(X, JY).$$

\subsection{Deformations and contractions of symplectic Lie algebras}

Let $(\g,\theta)$ be a symplectic Lie algebra. Recall that the symplectic form satisfies
\begin{enumerate}
  \item an open condition : $\theta^p \neq 0$ where $2p=\dim \g$,
  \item a closed condition : $d\theta =0$.
\end{enumerate}
The notion of formal deformation (see \cite{Ge,GRValued} ) is not adapted to these algebraic conditions. In fact, the $2p$-dimensional abelian Lie algebra $\mathcal{A}_{2p}$ is symplectic (any form is closed). But any $2p$-dimensional Lie algebra is isomorphic to a deformation of $\mathcal{A}_{2p}$ and we know many examples of $2p$-dimensional Lie algebras which are not symplectic.

The notion of contraction is better adapted. Recall the definition. Let $V$ be a finite dimensional $\K$-vector space and $f: ]0,1] \ra GL(V)$ be a continuous function. Let $[\ ,\ ]$ be a Lie bracket
on $V$. We defined the one parametrized family of  Lie brackets on $V$ by
\[
[x,y]_{\varepsilon}=f(\varepsilon)([f(\varepsilon)^{-1}(x), f(\varepsilon)^{-1}(y)]).
\]
The Lie bracket $[\ ,\ ]_\varepsilon$ is isomorphic to $[\ ,\ ]$.
\begin{definition}
If the limit
\[
\lbrack x,y\rbrack ^0 =\lim_{\varepsilon \rightarrow 0} [x,y]_{\varepsilon}
\]
exists, then $\lbrack , \rbrack ^0$ is a Lie bracket on $V$ and the Lie algebra
$(V,\lbrack , \rbrack ^0)$ is called a {\it contraction} of the Lie algebra $(V,[,])$.
\end{definition}

If we denote by $\mathcal{L}_n$ the algebraic variety whose elements are the Lie brackets on $\K^n$ (with the identification of a Lie bracket $\mu$ with its structure constants $\{ C_{ij}^k \}$   related with a fixed basis $\{e_1,\cdots,e_n\}$ of $\K^n$ by $\mu(e_i,e_j)=\sum C_{ij}^ke_k$ satisfying Jacobi polynomial conditions), then $\mathcal{L}_n$ is fibered by the orbits associated with the action of $GL(\K^n)$ on $\mathcal{L}_n$, the orbit of a Lie bracket $\mu$ is constituted of all  Lie bracket isomorphic to $\mu$. We denote by  $\mathcal{O}(\mu)$ the  orbit of $\mu$ in $\mathcal{L}_n$. We provide this algebraic variety $\mathcal{L}_n$  with the Zariski topology and we denote by $\overline{\mathcal{O}(\mu)}$ the closure of the orbit $\mathcal{O}(\mu)$. Let us note that when $\K=\C$, then $\overline{\mathcal{O}(\mu)}$ is also a differentiable manifold and the Zariski closure $\overline{\mathcal{O}(\mu)}$ in $\mathcal{L}_n$ coincides with the closure of  ${\mathcal{O}(\mu)}$ in the vector space corresponding to the set of constant structures $\{C_{ij}^k\}$ with the classical metric topology.

\begin{definition}
Let $\g=(\K^n, \mu)$ a $n$-dimensional $\K$-Lie algebra. A $n$-dimensional $\K$-Lie algebra $\g_0=(\K^n, \mu_0)$
is a degenrescence of $\g$ if
$\mu_0 \in \overline{\mathcal{O}(\mu)}$.
\end{definition}

It is easy to see that
any Lie algebra contraction of $(\K^n,\g)$ is a Lie algebra degeneration of this Lie algebra.

\medskip

Assume now that $\g=(\K^{2p}, \theta)$ is a symplectic Lie algebra. 
Let $\g_0=(\K^{2p},\mu_0)$ be a degenerescence of $\g$.
In general $\g _0$ is not symplectic. For example let $\g$ be the $6$-dimensional nilpotent Lie algebra given by
$$\left\{\begin{array}{l}
[ X_1,X_2] =X_3, \, \left[ X_1,X_3 \right]=X_4,
\left[ X_1,X_4 \right]=X_5,\left[ X_1,X_5 \right]=X_6, \\
\left[ X_2,X_3\right]=X_5,\left[ X_2,X_4 \right]=X_6,
\end{array}
\right.
$$
Denote by $\{ \alpha_i \}_i$ the dual basis. The form
 $\theta= \alpha_1 \wedge \alpha_6 +2 \alpha_2 \wedge \alpha_5 -\alpha_3 \wedge \alpha_4$ is symplectic.
Let $\{ Y_i\}_i$ be the new basis
$$\left\{
\begin{array}{l}
Y_i=\varepsilon X_i \ {\rm for} \ 1\leq i \leq 5, \\
Y_6=\varepsilon^2 X_6.
\end{array}
\right.
$$
%$$
%If $\{ \beta_i\}_i $ is the dual basis, that is,
%$$\left\{
%\begin{array}{l}
%\beta_i=\varepsilon^{-1} \alpha_i \ {\rm for} \ 1\leq i \leq 5 \\
%\beta_6=\varepsilon^{-2} \apha_6
%\end{array}
%\right.
%$$
We have
$$\left\{\begin{array}{l}
[ Y_1,Y_i] =\varepsilon^{-1}Y_{i+1}, \ i=2,3,4 \ \left[ Y_2,Y_3 \right]=\varepsilon^{-1}Y_5,\\
\left[ Y_1,Y_5 \right]=Y_6,\left[ Y_2,Y_4 \right]=Y_6.
\end{array}
\right.
$$
If $\varepsilon \rightarrow 0$ we obtain the contraction $\g_0$ of $\g$ defined by
$$\left\{
\begin{array}{l}
\left[ Y_1,Y_5 \right]=Y_6,\left[ Y_2,Y_4 \right]=Y_6
\end{array}
\right.
$$
which is isomorphic to $\h_5 \oplus \K$ where $\h_5$ is the $5$-dimensional Heisenberg algebra. We will see later that this Lie algebra is not symplectic.

However, we have

\begin{proposition}Let $(\g, \theta)$ be a symplectic Lie algebra and $\g_0$ be a contraction of $\g$ associated with the isomorphism $f_\varepsilon.$ If there exists $k$ such that $f_\varepsilon^*\theta=\varepsilon^k\theta$, then $(\g_0, \theta)$ is also symplectic.
\end{proposition}

\subsection{Symplectic nilpotent Lie algebras}

Since we are interested in this paper to the study of symplectic structures on some nilpotent Lie algebras,
we will recall briefly some basis facts concerning these Lie algebras.
The lower central series of a Lie algebra $\g$  is defined by
$$\rmc^0(\g)=\g,\quad  \rmc^j(\g)=[\g,\rmc^{j-1}(\g)],\;\; j\geq 1.
$$
Then $g$ is nilpotent if and only if  there exists $k$ with
$$\rmc^k(\mgg)=0.$$
The smallest index $k$ such that $\rmc^k(\mgg)=0$ and $\rmc^{k-1}(\mgg)\neq 0$ is the nilindex of $\g$. In this case we will say that $\g$ is $k$-step nilpotent. An interesting class of nilpotent Lie algebras is the class of $2$-step nilpotent Lie algebras, sometimes called metabelian Lie algebras (see, for example, \cite{Eberlein}).  For a nilpotent Lie algebra $\g$, we define its {\it characteristic sequence} which is an invariant up to isomorphism and that we recall in the following:
Let $X \neq 0$ be in $\g$. The linear operator $ad(X)$ is nilpotent and we denote by  $c(X)$ the decreasing sequence of the dimensions of Jordan blocks of $ad(X).$ The characteristic sequence of $\g$ is
 $$c(\g)= \max\{ c(X), X \in \g-\rmc^1(\g) \},$$
 considering the lexicographic order.
If $\g$ is $k$-step nilpotent, then  $c(\g)=(k,\cdots, k,k-1, \cdots, 1)$ where $k$ appears $n_k$ times, $k-1$ appears $n_{k-1}$ times, $\cdots$, and $1$ appears $n_1$ times with $n_k \neq 0$ and $n=\sum_{i=1}^{k} in_i.$ Since $ad(X)(X)=0$, we have also $n_1 \geq 1.$ For example, if $\g$ is $2$-step nilpotent, then $c(\g)=(2,\cdots,2,1,\cdots,1)$.

\medskip

The upper central series of $\mgg$  is defined by
$$\rmc_0(\mgg)=0,\quad  \rmc_j(\mgg)=\{X\in\mg:\,[X,\mg]\subset \rmc_{j-1}(\mgg)\},\;\; j\geq 1.
$$

\begin{theorem}\cite{dB3}

Let $\g=\rmc^0(\g)\supset\rmc^1(\g)\supset\dots$ and
$0=\rmc_0(\g)\subset\rmc_1(\g)\subset\dots$ be, respectively, the
upper and lower central series of a Lie algebra $\g$.
If  $\g$  is a symplectic Lie algebra then
\[
\dim\rmc_j(\g)+\dim \rmc^j(\g)\leq \dim \g
\]
for all $j\ge0$.
\end{theorem}

\medskip

We remark that this result doesn't give obstruction for $2$-step nilpotent Lie algebras with maximal  characteristic sequence that is equal to $(2,\cdots,2,1,1)$. In fact, from \cite{GRKegel} we have in this case  $\dim\rmc_j(\mgg)+\dim \rmc^j(\mgg)=  \dim \mgg$.

\medskip
\begin{theorem}
A nilpotent Lie algebra is never frobeniusian.
\end{theorem}

\pf This is a direct consequence of the following lemma proved in \cite{GRKegel,GNote}.

\begin{lemma}
Let $\g$ be a nilpotent Lie algebra over $\K$. Then the Cartan class of any non null linear form $\alpha \in \g^*$ is odd.
\end{lemma}

Recall that the Cartan class of a linear form can be computed with
\begin{enumerate}
\item$ cl(\alpha)=2p$ if and only if $(d\alpha)^p \neq 0, \alpha \wedge  (d\alpha)^p =0.$
\item $ cl(\alpha)=2p+1$ if and only if $ \alpha \wedge  (d\alpha)^p \neq 0, (d\alpha)^{p+1} =0.$
\end{enumerate}

If $\theta=d\alpha$ is symplectic, its Cartan class is even, equal to the dimension of $\g$. This is impossible.

\medskip

There exists a process to construct all the nilpotent symplectic Lie algebras. It is often called the construction by double extension and it was described by
 Alberto Medina and Philippe\ Revoy (\cite{MR}). Let $(\mathfrak{
g,\theta )}$ be a symplectic Lie algebra. We have seen that the product  $XY$ defined by
\[
\theta (XY,Z)=-\theta (Y,[X,Z])
\]
is left-symmetric, that is, satisfies (\ref{ls}). Let $D$ be a derivation of this left-symmetric algebra. Then the bilinear map
 $f$ on $\mathfrak{g}$ given by
\[
f(X,Y)=\theta (D(X),Y)+\theta (X,D(Y))
\]
is a closed $2$-form on $\mathfrak{g}$. Then we can define a one-dimensional central extension $\g_1=\g \oplus \K\{e\}$ of $\g$
by
\[
[X+\lambda e,Y+\beta e]_{\g_1}=[X,Y]+f(X,Y)e.
\]
If $D^*$ denotes the adjoint of $D$ with respect to the nondegenerate bilinear form $\theta$, that is
$$\theta(D(X),Y)=\theta(X,D^*(Y)),$$ then the bilinear form $g$ defined by
\[
g (X,Y)=\theta (((D+D^{\ast })\circ D+D^{\ast }\circ (D+D^{\ast
}))X,Y)
\]
is also a closed bilinear form. Assume that a closed form is an exact form, that is there exists $\alpha \in \g^*$ with $g(X,Y)=d\alpha (X,Y)= -\alpha [X,Y].$  There exists $X_g \in \g$ such that
\[
g(X,Y)=-\alpha [X,Y]= \theta (X_{g },[X,Y])
\]
for all $X,Y\in \mathfrak{g.}$  We can now define a new derivation $D_1$ of the Lie algebra $\g_1$ by
$$
\left\{
\begin{array}{l}
D_{1}(X) =-D(X)-\theta (X_{g},X)e, \ \forall X\in \mathfrak{g}, \\
D_{1}(e)=0,
\end{array}
\right.
$$
and we consider the one-dimensional extension by derivation $\g_2=\g_1 \oplus \K\{d\}$  of $\g_1$. Its bracket is given by
$$
\left\{
\begin{array}{l}
\lbrack X,Y]_{\g_2} =[X,Y]_{\g_1},\ \ \forall X,Y \in \g_1 \\
\lbrack d,X]_{\g_2} =-D_{1}(X)-\theta (X_{g },X)e
\end{array}
\right.
$$
We remark that $\dim \G_2= \dim \g+2.$ The bilinear map on $\g_2$
given by
\[
\left\{
\begin{array}{l}
\medskip
\theta _{1}\mid _{\mathfrak{g\times g}}=\theta \\
\theta _{1}(e,d)=1%
\end{array}%
\right.
\]%
the other non defined products being equal to $0$, is a symplectic form on $\g_2$. This last symplectic Lie algebra is called a double extension of the symplectic Lie algebra $\g$. Medina and Revoy have shown that any nilpotent symplectic Lie algebra is a double extension of a symplectic nilpotent Lie algebra.

%%%%%%%%%%%%%%%%%%%%%%%%%%%%%%%%%%%%%%%%%%%%%%%%%%%%%%%%
%%%%%%%%%%%%%%%%%%%%%%%%%%%%%%%%%%%%%%%%%%%%%%%%%%%%%%%
%%%%%%%%%%%%%%%%%%%%%%%%%%%%%%%%%%%%%%%%%%%%%%%%%%%%%%%
%%%%%%%%%%%%%%%%%%%%%%%%%%%%%%%%%%%%%%%%%%%%%%%%%%%%%%%

\section{$2$-step nilpotent Lie algebras of dimension $8$}

If $\g$ is a non abelian $8$-dimensional $2$-step nilpotent Lie algebra, its characteristic sequence is  one of the following $(2,2,2,1,1), (2,2,1,1,1,1), (2,1,1,1,1,1,1).$  If $\mathcal{F}^{8,2}$ is the set of such Lie algebras, we have defined in \cite{GRKegel}  a notion of restricted deformations that is deformations of Lie algebras of $\mathcal{F}^{8,2}$ remaining in this set. For this type of deformations, we have described the cohomology of deformations, denoted by $H_{CH}^*(\g,\g)$ for $\g \in \mathcal{F}^{8,2}.$ 

\subsection{Lie algebras of $\mathcal{F}^{8,2}_{3,2}$}
Let $\mathcal{F}^{8,2}_{3,2}$ be the subset of $\mathcal{F}^{8,2}$ corresponding to Lie algebras of characteristic sequence $(2,2,2,1,1)$.  Recall the description of this set  ( \cite{GRKegel}).     Let us denote by $\g_{3,2}$ the $8$-dimensional Lie algebra given by the brackets
$$
[X_1,X_{2i}]=X_{2i+1}, \quad 1\leq i \leq 3,
$$
and the other brackets are equal to zero or deduced by skew-symmetry.
Any $2$-step nilpotent $8$-dimensional Lie algebra with characteristic sequence $(2,2,2,1,1)$  is isomorphic to a linear deformation of $\g_{3,2}$. Its Lie bracket is isomorphic to $\mu=\mu_0 + t \varphi$,  where $\mu_0$ is the Lie bracket of $\g_{3,2}$, and $\varphi$ is a skew-bilinear map such that
\begin{equation}
\label{def}
\left\{
\begin{aligned}
&    \varphi \in Z^2_{CH}(\g_{3,2},\g_{3,2}), \\
&     \varphi \circ_1 \varphi=0
\end{aligned}
\right.
\end{equation}
where $(\varphi \circ_1 \varphi )(X,Y,Z)=\varphi ( \varphi (X,Y),Z).$
We deduce that
the family $\mathcal{F}^{8,2}_{3,2}$  of $2$-step nilpotent $8$-dimensional Lie algebras with characteristic sequence $(2,2,2,1,1)$
 is the union of two algebraic components, the first one, $\mathcal{C}_1(\mathcal{F}^{8,2}_{3,2})$, corresponding to the cocycles
\begin{equation}
\label{C1}
\varphi(X_{2i},X_{2j})= \sum_{k=1}^{3}a_{2i,2j}^{2k+1}X_{2k+1}, \quad 1\leq i<j \leq 4,
\end{equation}
the second one, $\mathcal{C}_2(\mathcal{F}^{8,2}_{3,2})$,  to the cocyles
\begin{equation}
\label{C_2}
\varphi(X_{2i},X_{2j})= \sum_{k=1}^{3}a_{2i,2j}^{2k+1}X_{2k+1}+a_{2i,2j}^{8}X_{8}, \quad 1\leq i<j \leq 3,
\end{equation}
where the undefined products $\varphi(X,Y)$ are nul or obtained by skew-symmetry.
Each of these components is a  regular algebraic variety. These  components can be characterized by the property of the number of generators. We have also mentioned that
 the Lie algebra $\mathfrak{h}_{8,1}$ of $\mathcal{C}_1(\mathcal{F}^{8,2}_{3,2})$ given by
  $$
  [X_1,X_{2i}]=X_{2i+1}, \;\;i=1,2,3, \quad [X_{2},X_{4}]=X_{7},  \quad [X_{4},X_{8}]=X_3,  \quad [X_{6},X_{8}]=X_5.
  $$
is rigid in the algebraic variety $2Nilp_{8}$ of $8$-dimensional $2$-step nilpotent Lie algebras.  Likewise,
 the Lie algebra $\mathfrak{h}_{8,2}$ of $\mathcal{C}_2(\mathcal{F}^{8,2}_{3,2})$ defined by
  $$
  [X_1,X_{2i}]=X_{2i+1}, \;\; i=1,2,3, \quad [X_{2},X_{6}]=X_{5},  \quad [X_{2},X_{4}]=X_8.
  $$
 is rigid in $2Nilp_{8}$.
  
\subsection{Lie algebras of $\mathcal{F}^{8,2}_{2,4}$}
Let $\mathcal{F}^{8,2}_{2,4}$ be the subset of $\mathcal{F}^{8,2}$ corresponding to Lie algebras of characteristic sequence $(2,2,1,1,1,1)$.  This case is not treated in \cite{GRKegel}. However, by proceeding in a similar way, we have immediately
that if  $\g_{2,4}$ the $8$-dimensional Lie algebra given by the brackets
$$
[X_1,X_{2i}]=X_{2i+1}, \quad i=1,2,
$$
then $2$-step nilpotent $8$-dimensional Lie algebra with characteristic sequence $(2,2,1,1,1,1)$  is isomorphic to a linear deformation of $\g_{2,4}$. Its Lie bracket is isomorphic to $\mu=\mu_0 + t \varphi$,  where $\mu_0$ is the Lie bracket of $\g_{2,4}$, and $\varphi$ is a skew-bilinear map such that
$$
\left\{
\begin{aligned}
&    \varphi \in Z^2_{CH}(\g_{2,4},\g_{2,4}), \\
&     \varphi \circ_1 \varphi=0.
\end{aligned}
\right.
$$

We deduce that
the family $\mathcal{F}^{8,2}_{2,4}$  of $2$-step nilpotent $8$-dimensional Lie algebras with characteristic sequence $(2,2,1,1,1,1)$
 is the union of two algebraic components, the first one, $\mathcal{C}_1(\mathcal{F}^{8,2}_{2,4})$, corresponding to the cocycles
\begin{equation}
\label{C3}
\varphi(X_{i},X_{j})= \sum_{k=1}^{2}a_{i,j}^{2k+1}X_{2k+1}, \quad 2\leq i<j \leq 8, \quad i,j \neq 3,5
\end{equation}
the second one, $\mathcal{C}_2(\mathcal{F}^{8,2}_{2,4})$,  to the cocyles
\begin{equation}
\label{C_2}
\left\{
\begin{array}{l}
\varphi(X_{i},X_{j})= \sum_{k=1}^{2}a_{i,j}^{2k+1}X_{2k+1},  \quad (i,j)\in \{(2,6),(2,7),(4,6),(4,7),(6,7)\}, \\
\varphi(X_{2},X_{4})=X_8
\end{array}
\right.
\end{equation}
where the undefined products $\varphi(X,Y)$ are nul or obtained by skew-symmetry.
Each of these components is a  regular algebraic variety. These  components can be characterized by the property of the number of generators.

\subsection{Lie algebras of $\mathcal{F}^{8,2}_{1,6}$}
Let $\mathcal{F}^{8,2}_{1,6}$ be the subset of $\mathcal{F}^{8,2}$ corresponding to Lie algebras of characteristic sequence $(2,1,1,1,1,1,1)$.  We know, from \cite{GRcontact} that any Lie algebras belonging to $\mathcal{F}^{8,2}_{1,6}$ is isomorphic to one of the following
\begin{itemize}
\item $\h_3 \oplus \A_5$
\item $\h_5 \oplus \A_3$
\item $\h_7 \oplus \A$
\end{itemize}
where $\h_{2p+1}$ is the $(2p+1)$-dimensional Heisenberg algebra and $\A_r$ the abelian ideal of dimension $r$.

%%%%%%%%%%%%%%%%%%%%%%%%%%%%%%%%%%%%%%%%%%%%%%%%%%%%%%%%
%%%%%%%%%%%%%%%%%%%%%%%%%%%%%%%%%%%%%%%%%%%%%%%%%%%%%%%%%
%%%%%%%%%%%%%%%%%%%%%%%%%%%%%%%%%%%%%%%%%%%%%%%%%%%%%%%

\section{Symplectic structures on $2$-step nilpotent Lie algebras}

\subsection{Recall some results}

We know many examples of  $2$-step nilpotent symplectic Lie algebras. The most recent \cite{dB1} based
on the  determination  of all nilradicals of
parabolic subalgebras of split real simple Lie algebras admitting symplectic structures are described into the following list.

\begin{enumerate}
 \item The abelian Lie algebras of dimesion even.
 \item $\K \oplus \n_{2,2}$ and $\n_{3,2}$ where $\n_{k,2}$ is the free $2$-step nilpotent Lie algebra generated by $k$ elements.
 \item $\{X_i,Y_i,\,Z_{ij}:1\leq i ,  j\leq 2\}$, $[X_i,Y_j]=Z_{ij}$, $1\!\leq\! i ,  j\!\leq\! 2$.
 \item $\{X_i,Y_i,\,Z_{ij}:1\!\leq\! i\!\le\! j\!\leq\! n\}$, $[X_i,Y_j]=[X_j,Y_i]=Z_{ij}$, $1\!\leq\! i\!\le\! j\!\leq\! n$; $n=2,3$.
 \item $\{X_i,Y_i,\,Z_{ij}:1\!\leq\! i\! <\!  j\!\leq\! 3\}$, $[X_i,Y_j]=-[X_j,Y_i]=Z_{ij}$, $1\!\leq\! i\! <\!  j\!\leq\! 3$.
 \item $\R X\ltimes_{\ad_X}(\R^{n-1}\oplus\R^{n-1})$
with $
\ad_X=\left(\begin{smallmatrix}
0&0\\[1mm]
I&0
\end{smallmatrix}\right)
$.
\end{enumerate}

Another interesting approach to determine symplectic structures on $2$-step nilpotent Lie algebra is presented in \cite{PO-TI}: given a finite graph $\mathcal{G}$, we consider the $\K$-vector space $V$ generated by the vertices and the $\K$-vector space generated by the edges. Let $\g_\mathcal{G}$ be the Lie algebra whose underlying vector space is $V_0 \oplus V_1$ and the Lie bracket defined by $[V_0,V_1]=[V_1,V_1]=0$ and $[X,Y]$ is the vector of $V_1$ if there exists an edge between the vertices $X$ and $Y$. This bracket satisfies the Jacobi identity and $\g_\mathcal{G}$ is a $2$-step nilpotent Lie algebra.
\begin{proposition} \cite{PO-TI}
Assume that $\g_\mathcal{G}$ is of dimension even. Then this Lie algebra is symplectic if and only if in each connected component of $\g_\mathcal{G}$, the number of edges doesn't exceed the number of vertices.
\end{proposition}

For example, we consider $\frak{n}_{k,2}$ the free $2$-step nilpotent Lie algebra with $k$ generators. We have a natural basis $\{X_1,\cdots,X_k,Y_{ij}=[X_i,X_j], 1 \leq i < j \leq k\}$.
There is a graph $\mathcal{G}$ such that $\frak{n}_{k,2}=\g_\mathcal{G}$. It is the polygon whose vertices are the generators $\{X_i\}$ and the edges all the diagonals. Then the numbers of vertices is $k$ and the number of edges $\frac{k(k-1)}{2}.$ From the previous proposition, $\frak{n}_{k,2}$ admits a symplectic form if and only if $k \leq 3$. Since $\frak{n}_{3,2}$ is of even dimension, it is the unique $2$-step nilpotent free Lie algebra which admits a symplectic form.

\medskip

To end this short list of known results, since the classification of $6$-dimensional symplectic nilpotent Lie algebra is established (see \cite{BG,GKNilp}), we deduce that

\begin{proposition}
Any $2$-step nilpotent Lie algebra of dimension $4$ is symplectic. Any $2$-step nilpotent Lie algebra of dimension $6$ which is not isomorphic to the direct product $\h_5 \oplus \A$ is symplectic where   $\h_5$ is the $5$-dimensional Heisenberg algebra.
\end{proposition}
In dimension $4$, up to an isomorphism we have only one Lie algebra which is the direct product of the $3$-dimensional Heisenberg algebra $\h_3$ by $\K$, that is with bracket
$$[X_1,X_2]=X_3$$
the other non defined being zero. In dimension $6$, the indecomposable are isomorphic to one of the following
\begin{enumerate}
  \item $[X_1, X_2] = X_3$, $[X_1, X_4] = X_5$, $[X_2, X_6] = X_5$ and  $\theta= \alpha_1\alpha_5+\alpha_2\alpha_4+\alpha_3\alpha_6$ is symplectic.
  \item  $[X_1, X_2] = X_3$, $[X_1, X_4] = X_5$, $[X_2, X_6] = X_5$, $[X_4,X_6]=X_3$ and  $\theta= \alpha_1\alpha_5+\alpha_2\alpha_4+\alpha_3\alpha_6$ is symplectic.
  \item If $\K=\R$,   $[X_1, X_2] = X_3$, $[X_1, X_4] = X_5$, $[X_2, X_6] = -X_5$, $[X_4,X_6]=X_3$ and  $\theta= \alpha_4\alpha_5+\alpha_1\alpha_6+\alpha_3\alpha_2$ is symplectic.

\end{enumerate}
The decomposable $2$-step nilpotent $\A_6$,  $\h_3\oplus \h_3$, $\h_3 \oplus \A_3$ are also symplectic where $\A_n$ is the $n$-dimensional abelian Lie algebra. But $\h_5 \oplus \A$ is not symplectic. In fact, $\h_5$ admits a linear contact form and its Reeb vector is in the Kernel of any closed skew-bilinear form of $\h_5 \oplus \A$.

\subsection{Case of dimension $8$ and characteristic sequence $(2,2,2,1,1)$}

\subsubsection{Case $n=8$ and $\g \in \mathcal{C}_2(\mathcal{F}^{8,2}_{3,2})$}

 The Lie bracket of a Lie algebra in this family has the general form
\begin{equation}\label{C1}
\left\{
\begin {array}{l}
[X_1,X_{2i}]=X_{2i+1}, i=1,2,3 \\
\lbrack X_{2i},X_{2j} \rbrack = \sum_{k=1}^3 C_{2i,2j}^{2k+1}X_{2k+1}+C_{2i,2j}^{8}X_{8}, \ 1\leq i<j \leq 3.
\end{array}
\right.
\end{equation}
If $C_{2i,2j}^{8} \neq 0$ for some pair $(i,j)$, then, up to an isomorphism, we can assume that
$$[X_2,X_4]=X_8 \ {\rm and } \ \  C_{2,6}^8=0=C_{4,6}^8.$$
We denote by $\{\alpha_1,\cdots,\alpha_8\}$ the dual basis of $\{X_1,\cdots,X_8\}$.  Let $\theta$ be a symplectic form on $\g$. The condition $d\theta =0$ implies
$$\theta(X_3,X_5)=\theta(X_3,X_7)=\theta(X_3,X_8)=\theta(X_5,X_7)=\theta(X_5,X_8)=\theta(X_7,X_8)=0$$
thus
$$\theta=\alpha_3 \wedge \beta_3 +\alpha_5 \wedge \beta_5+\alpha_7 \wedge \beta_7+\alpha_8 \wedge \beta_8$$
with $\beta_i= a_i^1 \alpha_1+a_i^2\alpha_2+a_i^4\alpha_4+a_i^6 \alpha_6$ for $i=3,5,7,8$.

\noindent We consider the $2$-form
$$\theta=\alpha_3 \wedge \alpha_6 +\alpha_5 \wedge ( \alpha_2 +\alpha_6) +  \alpha_7\wedge (\alpha_2+\alpha_4) +\alpha_8 \wedge (\alpha_1 + a_8^6\alpha_6) $$
with $a_8^6=C^7_{2,6}-C^5_{4,6}-C^7_{4,6}$. This form is symplectic on $\g$.

\begin{proposition}
Any $8$-dimensional  $2$-step nilpotent Lie algebra whose characteristic sequence is $(2,2,2,1,1)$ with $4$ generators (that is belonging to  $\mathcal{C}_2(\mathcal{F}^{8,2}_{3,2})$) is symplectic.
\end{proposition}

\subsubsection{Case $n=8$ and $\g \in \mathcal{C}_1(\mathcal{F}^{8,2}_{3,2})$}

Recall that the set  $\mathcal{C}_1(\mathcal{F}^{2p,2}_{p-1,2})$  is the family of $2$-step nilpotent Lie algebras isomorphic to Lie algebras given by
$$
\left\{
\begin{array}{l}
[X_1,X_{2i}]=X_{2i+1}, \ i=1,\cdots p-1,  \\
\lbrack X_{2i},X_{2j}\rbrack =\displaystyle \sum_{k=1}^{p-1}C_{2i,2j}^{2k+1}X_{2k+1}, \ 1\leq i<j \leq p.
\end{array}
\right.
$$
In particular $\mathcal{C}_1(\mathcal{F}^{8,2}_{3,2})$ is the closure of the rigid Lie algebra  $\h_{8,2}$ given by
$$[X_1,X_{2i}]=X_{2i+1}, i=1,2,3, \ [X_{2},X_{4}]=X_{7},  \ [X_{4},X_{8}]=X_3,  \ [X_{6},X_{8}]=X_5.$$

\begin{proposition}
The $8$-dimensional $2$-step nilpotent Lie algebra $\h_{8,2}$ which is rigid in $2Nilp_{8}$ is symplectic.
\end{proposition}
\pf Let $\theta= \alpha_1\wedge \alpha_8+\alpha_2\wedge \alpha_7-\alpha_3\wedge  \alpha_6+\alpha_4\wedge\alpha_5$ where the $\alpha_i$ constitue the dual basis of the basis $\{X_i\}$.
We have $d\theta= 0$ and $\theta ^4 \neq 0$. The form $\theta$ is a symplectic.  \qed

\medskip

Let $g \in \mathcal{C}_1(\mathcal{F}^{8,2}_{3,2})$. It is a contraction of $\h_{8,2}$. If the contraction is associated with the sequence of isomorphisms $f_\varepsilon$ such that $\lim f_\varepsilon^* \theta$  is always of rank $4$, then $\g$ is symplectic. Now, we have to determinate all the other cases. If $\g$ is symplectic, any symplectic form $\theta $ on $\g$ satisfies
$$\theta(X_i,X_j)=0$$
for $i,j=3,5,7.$  Consider the decomposition  $\theta=\theta_1+\theta_2$ with 
$\theta_1=\alpha_3 \wedge \beta_3 + \alpha_5 \wedge \beta_5 + \alpha_7 \wedge \beta_7$
and $\beta_i(X_j)=0$ for $i,j=3,5,7$ and $\theta_2=\sum a_{i,j} \alpha_i \wedge \alpha_j$ with $i,j \in I=\{1,2,4,6,8\}$. Since $d\alpha_i=0$ for any $i \in I$, we have $d\theta_2=0$. Then $\theta$ symplectic implies $d\theta_1=0$.
This implies that, in order to determinate symplectic Lie algebras, we have only to considerate on $\g$ skew bilinear forms $\theta_1$ with
$$\theta_1=\alpha_3 \wedge \beta_3 + \alpha_5 \wedge \beta_5 + \alpha_7 \wedge \beta_7$$ such that
 $d\theta_1=0$ and $ \beta_3\wedge \beta_5  \wedge \beta_7 \neq 0.$
 Since $d\beta_i=0$, then $d\theta_1=0$ is equivalent to
$$d\alpha_3 \wedge \beta_3 +d\alpha_5 \wedge \beta_5 +d\alpha_7 \wedge \beta_7=0.$$
So instead of  determining symplectic Lie algebras belonging to $\mathcal{C}_1(\mathcal{F}^{8,2}_{3,2})$, we will characterize the Lie algebras of this family such that any closed $2$-forms is not of maximal rank, that is satisfies $\theta^4=0$. This last condition is equivalent to $\theta_1^3 =0,$ that is, $ \beta_3\wedge \beta_5  \wedge \beta_7 =0$.  The linear forms $\beta_i, i=3,5,7$ are not independent and we can assume that $\beta_7=a\beta_3+b\beta_5$. Thus
$$\theta_1=(\alpha_3 +a\alpha_7)\wedge \beta_3 + (\alpha_5+b\alpha_7) \wedge \beta_5 $$
and there exists  a change of basis respecting the characteristic sequence giving 
$$\theta_1=\alpha_3 \wedge \beta_3 + \alpha_5 \wedge \beta_5 $$
that is $\beta_7=0$. Then $\g$ is not symplectic if the condition $d\theta_1 =0$ implies $\beta_7=0$. In a first step, we assume that
$\beta_i=\ds\sum_{j=2,4,6,8}a_{i,j}\alpha_j.$ 
The condition $d\theta_1=0$ implies that $\beta_i(X_8)=0$ for $i=3,5,7$, that is $a_{i,8}=0,$ and $\beta_i=\ds\sum_{j=2,4,6}a_{i,j}\alpha_j.$
Let $M_i$ be the matrix
$$
M_i=
\left(
\begin{array}{ccc}
C_{46}^i & -C_{26}^i & C_{24}^i  \\
C_{48}^i & -C_{28}^i & 0   \\
C_{68}^i & 0 & -C_{28}^i \\
0 & C_{68}^i &   -C_{48}^i \\
\end{array}
\right)
$$
for $i=3,5,7.$
The conditions $d\theta_1(X_{2i},X_{2j},X_{2k})=0$ are equivalent to
$$M_7A_7=-M_3A_3-M_5A_5$$
where
$A_i=^t(a_{i,2},a_{i,4},a_{i,6})$, $i=3,5,7.$ If we want that $(d\theta_1=0 \Rightarrow \beta_7=0)$ we have that $(d\theta_1=0 \Rightarrow A_7=0)$, then $-M_3A_3-M_5A_5=0$ and the rank of $M_7$ is $3$. So if $\g$ is not symplectic, it implies 
that one of the expressions $C^7_{2,8} \Delta$, $C^7_{4,8} \Delta$, $C^7_{6,8} \Delta$ is not zero, where $\Delta=C^7_{2,8}C^7_{46}+C^7_{2,4}C^7_{6,8}-C^7_{2,6}C^7_{4,8}$. The condition $\Delta \neq 0$ is equivalent to the say that the $2$-form $d\alpha_7-\alpha_1 \wedge \alpha_6$ is of (maximal) rank $4$.  If moreover $C_{6,8}^7  \neq 0$,  up to isomorphism adapted to the characteristic sequence, we can reduce the form $d\alpha_7$ to $$d\alpha_7=-\alpha_1 \wedge \alpha_6-\alpha_2 \wedge \alpha_4-\alpha_6 \wedge \alpha_8$$
%corresponding to $C^7_{6,8}\neq 0$
and $\g$ is isomorphic to the Lie algebra given by
$$
\left\{
\begin{array}{l}
d\alpha_3=-\alpha_1 \wedge \alpha_2, \\
d\alpha_5=-\alpha_1 \wedge \alpha_4, \\
d\alpha_7=-\alpha_1 \wedge \alpha_6- \alpha_2 \wedge \alpha_4-\alpha_6 \wedge \alpha_8\\
d\alpha_i=0, \ i=1,2,4,6,8.
\end{array}
\right.
$$
We verify directly that this Lie algebra doesn't admits general symplectic form (corresponding to $a_{i,1} \neq 0$ for some $i \in \{3,5,7\}$).
In all the other cases, there exists a symplectic form with $a_{7,1}=1$.  We deduce:

\begin{theorem}
A $8$-dimensional $2$-step nilpotent Lie algebra whose characteristic sequence  is $(2,2,2,1,1)$ is symplectic if and only if it is not isomorphic to the Lie algebra whose Maurer-Cartan equations are
\begin{equation}
\label{nosym}
\left\{
\begin{array}{l}
d\alpha_3=-\alpha_1 \wedge \alpha_2, \\
d\alpha_5=-\alpha_1 \wedge \alpha_4, \\
d\alpha_7=-\alpha_1 \wedge \alpha_6- \alpha_2 \wedge \alpha_4-\alpha_6 \wedge \alpha_8\\
d\alpha_i=0, \ i=1,2,4,6,8.
\end{array}
\right.
\end{equation} 
\end{theorem}

\medskip

%%%%%%%%%%%%%%%%%%%%%%%%%%%%%%%%%%%%%%%%%%%
\subsection{Case of dimension $8$ and characteristic sequence $(2,2,1,1,1,1)$}

Let $\mathcal{F}^{8,2}_{2,4}$ be the family of $8$-dimensional $2$-step nilpotent Lie algebras whose characteristic sequence is equal to $(2,2,1,1,1,1)$.  We have seen that this set is constituted of two components,  $\mathcal{C}_1(\mathcal{F}^{8,2}_{2,4})$ whose elements are Lie algebras with $6$ generators and $\mathcal{C}_2(\mathcal{F}^{8,2}_{2,4})$ whose elements are Lie algebras with $5$ generators. 

\subsubsection{$n=8$ and $\g \in \mathcal{C}_1(\mathcal{F}^{8,2}_{2,4})$} Let us consider the first component $\mathcal{C}_1(\mathcal{F}^{8,2}_{2,4})$. Let $\mathfrak{k}_8 \in \mathcal{C}_1(\mathcal{F}^{8,2}_{2,4})$ whose Lie bracket is defined by the corresponding Maurer-Cartan equations
\begin{equation}
\label{k8}
\left\{
\begin{array}{l}
d\alpha_3=\alpha_1\wedge\alpha_2+\alpha_4\wedge\alpha_6+\alpha_7\wedge\alpha_8    \\
d\alpha_5=\alpha_1\wedge\alpha_4+\alpha_2\wedge\alpha_8+\alpha_6\wedge\alpha_7,
\end{array}
\right.
\end{equation}
where $\{\alpha_1,\cdots,\alpha_8\}$ is the dual basis of $\{X_1,\cdots,X_8\}$.
\begin{proposition}
The Lie algebra $\mathfrak{k}_8$ is rigid in $2Nilp_8$.
\end{proposition}
\pf Let us note, in a first step, that $\mathfrak{k}_8$ cannot be deformed in $\mathcal{F}^{8,2}_{3,2}$, that is, in a Lie algebra with characteristic sequence $(2,2,2,1,1)$. In fact, in $\mathfrak{k}_8$, the linear forms $\alpha_3$ and $\alpha_5$ have a Cartan class equal to $7$ (see, for the notions of Cartan classes \cite{GRcontact}). But in any Lie algebra of  $\mathcal{F}^{8,2}_{3,2}$, the center is of dimension greater than or equal to 3, this implies that the Cartan class of any linear form at most $6$. Since the maximal class is growing in the deformation process, this show that any deformation in $2Nilp_8$ of $\mathfrak{k}_8$ remains in $\mathcal{F}^{8,2}_{2,4}$. Then, we have to study deformations of $\mathfrak{k}_8$ in $\mathcal{F}^{8,2}_{2,4}$. We consider the deformation cohomology $H^*_{CH}(\mk_8,\mk_8)$ introduced in \cite{GRKegel}.
We can assume that any $2$-cocycle $\varphi$ satisfies $\varphi(X_1,Y)=0$ for any $Y$. Then $\dim Z^2_{CH}(\mk_8,\mk_8)= 20$ because we have
$$\varphi(X_i,X_j)=a_{ij}^3X_3+a_{ij}^5X_5$$
for $i=2,4,6,7$ and $j=4,6,7,8$ and $i<j$. Let $f$ be a linear endomorphism of $\mk_8$. Assume that $\delta_{CH}f(X_1,X)=0.$ This implies that
$$f(X_3)=[f(X_1),X_2]+[X_1,f(X_2)], \ f(X_5)=[f(X_1),X_4]+[X_1,f(X_4)]$$
and
$$[f(X_1),X_i]+[X_1,f(X_i)]=0, \ i=3,5,6,7,8.$$
This last identity is satisfied for $i=3$ and $i=5$ because $X_3$, $X_5$, $f(X_3)$ and $f(X_5)$ are in the center. The other identities imply
$$b_2^6=-b_4^1, b_4^6=b_7^1,b_2^7=b_8^1,b_4^7=-b_6^1,b_2^8=-b_7^1,b_4^8=-b_2^1.$$
Now, if we compute $\delta_{CH} \varphi (X_i,X_j)$ for $i=2,4,6,7$ and $j=4,6,7,8$ and $j>i$, writing $ \varphi (X_i,X_j)=a_{ij}^3X_3+a_{ij}^5X_5$, we see that we have no relations between all these coefficients. We deduce that $\dim B^2_{CH}(\mk_8,\mk_8)= 20$ and $H^2_{CH}(\mk_8,\mk_8)=0.$ Then this Lie algebra is rigid in $2Nilp_8$.

\begin{corollary}
The component $\mathcal{C}_1(\mathcal{F}^{8,2}_{2,4})$ is the closure of the orbit in $2Nilp_8$ of $\mk_8$. Any Lie algebra of $\mathcal{C}_1(\mathcal{F}^{8,2}_{2,4})$ not isomorphic to $\mk_8$ is a contraction of $\mk_8$.
\end{corollary}

\begin{proposition}
The Lie algebra $\mk_8$ is not symplectic.
\end{proposition}
\pf Let $\theta$ be a skew-symmetric bilinear map on $\mk_8$. We put $\theta=\theta_1+\theta_2$ with
$$\theta_2= \ds\sum a_{ij}\alpha_i \wedge \alpha_j, \ \ {\text {\rm with\ }} 1\leq i < j \leq 8, {\rm and\ } i,j \neq 3,5$$
and
$$\theta_1=\alpha_3 \wedge(\sum_{i=1}^8 a_i \alpha_i)+\alpha_5 \wedge(\sum_{i=1}^8 b_i \alpha_i)$$
with $a_3=b_3=b_5=0.$
We have $d\theta_2=0$, so $\theta$ is closed if and only if $d\theta_1=0.$ But
$$d\theta_1=(\alpha_1\wedge\alpha_2+\alpha_4\wedge\alpha_6+\alpha_7\wedge\alpha_8 )\wedge(\sum_{i=1}^8 a_i \alpha_i)+(\alpha_1\wedge\alpha_4+\alpha_2\wedge\alpha_8+\alpha_6\wedge\alpha_7) \wedge(\sum_{i=1}^8 b_i \alpha_i)$$
and $d\theta_1=0$ if and only if $a_i=b_i=0$ for any $i=1,\cdots,8.$ Then $\theta_1=0,$ $\theta=\theta_2$. But $\theta_2$ is not of maximal rank because $X_3$ and $X_5$ are in the kernel of this form. This implies that $\theta$ which is closed, is not symplectic.

\medskip

\noindent{\bf Remark.} Some contractions of $\mk_8$ can be writen
$$
\left\{
\begin{array}{l}
d\alpha_3=\alpha_1\wedge\alpha_2+\rho_1\alpha_4\wedge\alpha_6+\rho_2\alpha_7\wedge\alpha_8    \\
d\alpha_5=\alpha_1\wedge\alpha_4+\rho_3\alpha_2\wedge\alpha_8+\rho_4\alpha_6\wedge\alpha_7,   \\
d\alpha_i=0, \ i\neq 3,5
\end{array}
\right.
$$
where the parameters $\rho_i$ are parameters of contraction.  By  simple computation, we can see that the symplectic Lie algebras of this family are  Lie algebras with
$$\rho_1=\rho_4=0, \ {\text{\rm or }} \ \rho_2=\rho_3=0, \ {\text{\rm or }} \ \rho_2=\rho_4=0.$$
In each one of these cases, the corresponding Lie algebra is decomposable. We will show that this property is true in any (even) dimension greater than $8$.

\medskip

\begin{proposition}
Any $8$-dimensional $2$-step nilpotent symplectic Lie algebra with characteristic sequence $(2,2,1,1,1,1)$ and belonging to $\mathcal{C}_1(\mathcal{F}^{8,2}_{2,4})$
is decomposable.
\end{proposition}

\pf
 We have seen that the general Maurer-Cartan system of Lie algebras of $\mathcal{C}_1(\mathcal{F}^{8,2}_{2,4})$ is
\begin{equation}
\label{MC4}
\left\{
\begin{array}{l}
d\alpha_3=\alpha_1\wedge\alpha_2+\displaystyle\sum_{\begin{array}{l} 2 \leq i<j\leq 8\\ i,j \neq 3,5
\end{array}}a^3_{ij} \alpha_i\wedge\alpha_j   \\
d\alpha_5=\alpha_1\wedge\alpha_4+\displaystyle\sum_{\begin{array}{l} 2 \leq i<j\leq 8\\ i,j \neq 3,5
\end{array}}a^5_{ij} \alpha_i\wedge\alpha_j   \\
d\alpha_i=0, \ i=1,2,4,6,7,8,
\end{array}
\right.
\end{equation}
Let $\g$ be a Lie algebra of $\mathcal{C}_1(\mathcal{F}^{8,2}_{2,4})$. Assume that $\g$ is symplectic and let $\theta$ be a symplectic form on $\g$. Since $d\theta=0$, we have in particular
$$d\theta(X_1,X_2,X_5)=\theta(X_3,X_5)=0.$$
Then, $\theta=\theta_1+\theta_2$ with 
$$
\left\{
\begin{array}{l}
\theta_1=\alpha_3\wedge\beta_3+\alpha_5\wedge\beta_5,\\
\medskip
\theta_2=\displaystyle\sum_{\begin{array}{l} 1 \leq i<j\leq 8\\ i,j \neq 3,5
\end{array}}\rho_{i,j} \alpha_i\wedge\alpha_j   \\
\end{array}
\right.
$$
with
$\beta_3=\displaystyle\sum_{\begin{array}{l} 1 \leq i<j\leq 8\\ i,j \neq 3,5
\end{array}}a_i \alpha_i$ and $\beta_5=\displaystyle\sum_{\begin{array}{l} 1 \leq i<j\leq 8\\ i,j \neq 3,5
\end{array}}b_i \alpha_i$. It is clear that $d\theta_2=0$, then $\theta$ is symplectic if and only if $d\theta_1=0$ and $\theta_1^2 \neq 0.$ We will study the forms $\theta_1$. 

$\bullet$ Assume that $a_1$ or $b_1 \neq 0$. We can consider that $a_1 \neq 0$ and more precisely $a_1=1.$ We have
$$
\left\{
\begin{array}{lll} 
\theta_1 &=&\alpha_3\wedge(\alpha_1+a_2\alpha_2+a_4\alpha_4+a_6\alpha_6+a_7\alpha_7+a_8\alpha_8)+
\alpha_5\wedge(b_1\alpha_1+b_2\alpha_2+b_4\alpha_4\\
&&+b_6\alpha_6+b_7\alpha_7+b_8\alpha_8)\\
&=& (\alpha_3+b_1\alpha_5)\wedge\alpha_1+\alpha_3\wedge(a_2\alpha_2+a_4\alpha_4+a_6\alpha_6+a_7\alpha_7+a_8\alpha_8)+
\alpha_5\wedge(b_2\alpha_2\\
&&+b_4\alpha_4+b_6\alpha_6+b_7\alpha_7+b_8\alpha_8)\\
&=& (\alpha_3+b_1\alpha_5)\wedge(\alpha_1+a_2\alpha_2+a_4\alpha_4+a_6\alpha_6+a_7\alpha_7+a_8\alpha_8)+
\alpha_5\wedge((b_2-b_1a_2)\alpha_2\\
&&+(b_4-b_1a_4)\alpha_4+(b_6-b_1a_6)\alpha_6+(b_7-b_1a_7)\alpha_7+(b_8-b_1a_8)\alpha_8)\\
\end{array}
\right.
$$
This shows that, considering the change of basis $\alpha'_2=\alpha_2+b_1\alpha_4,\alpha'_3=\alpha_3+b_1\alpha_5$ and $\alpha'_i=\alpha_i$ for all the other elements of the basis, we can assume that
$$\theta_1=\alpha_3\wedge\beta_3+\alpha_5\wedge\beta_5$$
with
$$a_1=1, \ b_1=0.$$
Since $a_1=1$, we consider $\widetilde{\alpha_1}=\beta_3$ and $\widetilde{\alpha_i}=\alpha_i$ for all the other indices. This change of basis doesn't affect the Maurer-Cartan system (\ref{MC4}). Then, we can assume, in the original basis, that
$$\theta_1=\alpha_3\wedge\alpha_1+\alpha_5\wedge\beta_5$$
with
$$\beta_5=b_2\alpha_2+b_4\alpha_4+b_6\alpha_6+b_7\alpha_7+b_8\alpha_8.$$
The condition $d\theta_1=0$ implies, in particular that 
$$a^3_{26}=a^3_{27}=a^3_{28}=a^3_{67}=a^3_{68}=a^3_{78}=0$$
that is
$$d\alpha_3=\alpha_1\wedge\alpha_2+a^3_{24}\alpha_2\wedge\alpha_4+a^3_{46}\alpha_4\wedge\alpha_6
+a^3_{47}\alpha_4\wedge\alpha_7+a^3_{48}\alpha_4\wedge\alpha_8.$$
If one of the parameters $b_6,b_7,b_8$ is not zero, for example $b_6$,  a change of basis which preserves the previous reduction gives
$$\beta_5=b_2\alpha_2+b_4\alpha_4+b_6\alpha_6$$
and putting $\widetilde{\alpha}_6=b_2\alpha_2+b_4\alpha_4+b_6\alpha_6$, the expression of $d\alpha_3$ is preserved. This shows that, in this case, if $\g$ is symplectic, there exists a symplectic form associated with
$$\theta_1=\alpha_3\wedge\alpha_1+\alpha_5\wedge \alpha_6.$$
In this case $d\theta_1=0$ is equivalent to $a_{24}^3=a^3_{47}=a^3_{48}=0 ,1=-a_{46}^3,a_{24}^5=a_{27}^5=a_{28}^5=a_{47}^5=a_{48}^5=a_{78}^5=0.$ 
$$
\left\{
\begin{array}{l}
d\alpha_3=\alpha_1 \wedge \alpha_2 - \alpha_4 \wedge \alpha_6,\\
d\alpha_5= \alpha_1 \wedge \alpha_4 + (a^5_{26}\alpha_2+a^5_{46}\alpha_4-a^5_{67}\alpha_7-a^5_{68}\alpha_8)\wedge \alpha_6.
\end{array}
\right.
$$

If $b_6=b_7=b_8=0$, that is $\beta_5=b_2\alpha_2+b_4\alpha_4$, then $\theta_1=\alpha_3\wedge\alpha_1+\alpha_5\wedge(b_2\alpha_2+b_4\alpha_4)$ and $d\theta_1=0$ implies 
$d\alpha_3= \alpha_1\wedge\alpha_2+a^3_{24}\alpha_2\wedge\alpha_4$ and $b_2=a^3_{24}$ and $\theta_1^2 \neq 0$ implies $a^5_{67}=a^5_{68}=a^5_{78}=0.$ This permits to consider $a^5_{67}=a^5_{68}=a^5_{78}=0.$ Then the linear space generated by ${X_6,X_7,X_8}$ is an abelian subalgebra of $\g$. The endomorphisms $adX_i$, $i=6,7,8$ acts non trivially only on the space generated by $X_2$ and $X_4$. Thus, we can consider that there exists $i \in \{6,7,8\}$ with $adX_i=0$ and the Lie algebra $\g$ is decomposable.

$\bullet$ Assume that $a_1=b_1= 0$. Similar computations show that $\theta_1=\alpha_3\wedge(a_2\alpha_2+a_4\alpha_4)+\alpha_5\wedge(b_2\alpha_2+b_4\alpha_4)$ with $a_2b_4-a_4b_2 \neq 0$ and the Maurer-Cartan equations can be reduced to one of the following system:
$$
\left\{
\begin{array}{l}
d\alpha_3=\alpha_1\wedge\alpha_2+\alpha_2\wedge\alpha_6
, \\
d\alpha_5=\alpha_1\wedge\alpha_4+(a^5_{24}\alpha_2-a^5_{46}\alpha_6-a^5_{47}\alpha_7-a^5_{48}\alpha_8)
\wedge\alpha_4\\
d\alpha_i=1,2,4,6,7,8.
\end{array}
\right.
$$ or
$$
\left\{
\begin{array}{l}
d\alpha_3=\alpha_1\wedge\alpha_2+\alpha_4\wedge\alpha_6
, \\
d\alpha_5=\alpha_1\wedge\alpha_4+\alpha_2\wedge (a^5_{24}\alpha_4+a^5_{26}\alpha_6+a^5_{27}
\alpha_7+a^5_{28}\alpha_8
)\\
d\alpha_i=1,2,4,6,7,8.
\end{array}
\right.
$$
In these two cases, $\g$ is decomposable. We deduce the proposition.

\subsubsection{$n=8$ and $\g \in \mathcal{C}_2(\mathcal{F}^{8,2}_{2,4})$}  The general Maurer-Cartan system of Lie algebras of $\mathcal{C}_2(\mathcal{F}^{8,2}_{2,4})$ is
$$
\left\{
\begin{array}{l}
d\alpha_3=\alpha_1\wedge\alpha_2+\displaystyle\sum_{\begin{array}{l} 2 \leq i<j\leq 7\\ i,j \neq 3,5
\end{array}}a^3_{ij} \alpha_i\wedge\alpha_j   \\
d\alpha_5=\alpha_1\wedge\alpha_4+\displaystyle\sum_{\begin{array}{l} 2 \leq i<j\leq 7\\ i,j \neq 3,5
\end{array}}a^5_{ij} \alpha_i\wedge\alpha_j   \\
d\alpha_8=\alpha_2\wedge\alpha_4 \\
d\alpha_i=0, \ i=1,2,4,6,7,
\end{array}
\right.
$$
If some constant $a_{ij}^3$ or $a_{ij}^5$ is not zero for $i$ or $j$ greater than or equal to $6$, then the characteristic sequence is not $(2,2,1,1,1,1)$. Then $\g$ is isomorphic to the Lie algebra
\begin{equation}
\label{MC41}
\left\{
\begin{array}{l}
d\alpha_3=\alpha_1\wedge\alpha_2, \\
d\alpha_5=\alpha_1\wedge\alpha_4,  \\
d\alpha_8=\alpha_2\wedge\alpha_4 \\
d\alpha_i=0, \ i=1,2,4,6,7,
\end{array}
\right.
\end{equation}
and $\g$ is decomposable.

\subsubsection{Conclusion}

\begin{theorem}
Any symplectic Lie algebra with characteristic sequence $(2,2,1,1,1,1)$ is decomposable.
\end{theorem}

\medskip

We deduce immediately the description of symplectic Lie algebras with characteristic sequence $(2,2,1,1,1,1)$. The notation corresponds to the list given in
\begin{center}
www.livres-mathematics.fr
\end{center}
\noindent and the notation $\g_1 \oplus \g_2$ corresponds to a direct product of two ideals.

 1. Complex case
 \begin{itemize}
  \item $\n_7^{124} \oplus \C$ with Maurer-Cartan equations
$$
\left\{
\begin{array}{l}
d\alpha_3=\alpha_1\wedge\alpha_2+\alpha_4\wedge\alpha_7, \\
d\alpha_5=\alpha_1\wedge\alpha_4+\alpha_6\wedge\alpha_7,  \\
d\alpha_i=0, \ i=1,2,4,6,7,8
\end{array}
\right.
$$
  \item $\n_6^{19} \oplus \C^2$ with Maurer-Cartan equations
$$
\left\{
\begin{array}{l}
d\alpha_3=\alpha_1\wedge\alpha_2, \\
d\alpha_5=\alpha_1\wedge\alpha_4+\alpha_2\wedge\alpha_6,  \\
d\alpha_i=0, \ i=1,2,4,6,7,8
\end{array}
\right.
$$
    \item $\n_6^{20} \oplus \C^2$ with Maurer-Cartan equations
$$
\left\{
\begin{array}{l}
d\alpha_3=\alpha_1\wedge\alpha_2, \\
d\alpha_5=\alpha_1\wedge\alpha_4,  \\
d\alpha_6=\alpha_2\wedge\alpha_4,  \\
d\alpha_i=0, \ i=1,2,4,7,8
\end{array}
\right.
$$
 \item $\n_5^{5} \oplus \C^3$ with Maurer-Cartan equations
$$
\left\{
\begin{array}{l}
d\alpha_3=\alpha_1\wedge\alpha_2, \\
d\alpha_5=\alpha_1\wedge\alpha_4,  \\
d\alpha_i=0, \ i=1,2,4,6,7,8
\end{array}
\right.
$$
 \item $\n_3^{1}\oplus \n_3^1 \oplus \C^3$ with Maurer-Cartan equations
$$
\left\{
\begin{array}{l}
d\alpha_3=\alpha_1\wedge\alpha_2, \\
d\alpha_5=\alpha_1\wedge\alpha_4+\alpha_4\wedge\alpha_6,  \\
d\alpha_i=0, \ i=1,2,4,6,7,8.
\end{array}
\right.
$$
\end{itemize}

\medskip

 2. Real case.
\begin{itemize}
  \item $\n_7^{124} \oplus \R$ 
    \item $\n_6^{19} \oplus \R^2$
        \item $\n_6^{20} \oplus \R^2$
    \item $\n_6^{20,1} \oplus \R^2$ with Maurer-Cartan equations
$$
\left\{
\begin{array}{l}
d\alpha_3=\alpha_1\wedge\alpha_2+\alpha_4\wedge\alpha_6, \\
d\alpha_5=\alpha_1\wedge\alpha_4-\alpha_2\wedge\alpha_6,  \\
d\alpha_i=0, \ i=1,2,4,6,7,8
\end{array}
\right.
$$
 \item $\n_5^{5} \oplus \R^3$,
  \item $\n_3^{1}\oplus \n_3^1 \oplus \R^3$.
  \end{itemize}

\subsection{Case of dimension $8$ and characteristic sequence $(2,1,1,1,1,1,1)$}

A Lie algebra with characteristic sequence $(2,1,1,1,1,1,1)$ is isomorphic to $\h_3 \oplus \K^5$ or $\h_5 \oplus \K^3$ or $\h_7 \oplus \K$ (\cite{GKNilp}). We deduce

\begin{theorem}
Any $8$-dimensional $2$-step nilpotent Lie algebra with characteristic sequence $(2,1,1,1,1,1,1)$ is symplectic if and only if it is isomorphic to $\h_3 \oplus \K^5$.
\end{theorem}

%%%%%%%%%%%%%%%%%%%%%%%%%%%%%%%%%%%%%%%%%%%%%%%%%%%%%%%%%%%%%%%%%%%%%%
%%%%%%%%%%%%%%%%%%%%%%%%%%%%%%%%%%%%%%%%%%%%%%%%%%%%%%%%%%%%%%%%%%%%%%%
%%%%%%%%%%%%%%%%%%%%%%%%%%%%%%%%%%%%%%%%%%%%%%%%%%%%%%%%%%%%%%%%%%%%%%

\section{Leading to the general case}

\subsection{Symplectic Lie algebras in $\mathcal{F}^{2p,2}_{1,2p-2}$}
 Let $\g$ be a $2p$-dimensional symplectic nilpotent Lie algebras with characteristic sequence $(2,1,\cdots,1)$. We have recalled that such a Lie algebra is a direct sum of an Heisenberg algebra and an abelian ideal. From the previous section, we can affirm

\begin{theorem}
Any $2p$-dimensional $2$-step nilpotent Lie algebra with characteristic sequence $(2,1,\cdots,1)$ is symplectic if and only if it is isomorphic to $\h_3 \oplus \K^{2p-3}$. 
\end{theorem}

\subsection{Symplectic Lie algebras in $\mathcal{F}^{2p,2}_{2,2p-4}$}
Let $\g$ be a $2p$-dimensional  nilpotent Lie algebras with characteristic sequence $(2,2,1,\cdots,1)$.  We have seen that $\mathcal{F}^{2p,2}_{2,2p-4}$ is the union of two components, the first corresponding to Lie algebras with $2p-2$ generators and the second to Lie algebras with $2p-3$ generators. The second case is very simple to describe. From the previous section, we know that any Lie algebra belonging to $\mathcal{C}_2(\mathcal{F}^{2p,2}_{2,2p-4})$ is isomorphic to $\n_6^{20} \oplus \K^{2n-6}$ and these Lie algebras are symplectic. Assume now that $\g$ is in $\mathcal{C}_1(\mathcal{F}^{2p,2}_{2,2p-4})$. We have seen, using computational method, that for $p=4$, $\g$ is decomposable. We will show the same result, but for any  $p \geq 4$, using geometrical method.

The Maurer-Cartan equations of $\g$ can be reduced to
\begin{equation}
\label{k2p}
\left\{
\begin{array}{l}
d\alpha_3=\alpha_1\wedge\alpha_2+\sum a_{ij}^3\alpha_i\wedge\alpha_j, \ \ 2 \leq i <j \leq 2p, \ i,j \neq 3, 5  \\
d\alpha_5=\alpha_1\wedge\alpha_4+\sum a_{ij}^5\alpha_i\wedge\alpha_j, \ \ 2 \leq i <j \leq 2p, \ i,j \neq 3, 5,
d\alpha_i=0, \ i \neq 3,5.
\end{array}
\right.
\end{equation}
Let $\theta$ be a symplectic form on $\g$. It can be decomposed in $\theta=\theta_1+\theta_2$ with $\theta_1=\alpha_3 \wedge \beta_3+\alpha_5 \wedge \beta_5$, $d\theta_1=0$ and $\theta_1^2\neq 0$.  This implies
$$
\left\{
\begin{array}{l}
d\alpha_3\wedge \beta_3+d\alpha_5 \wedge \beta_5=0,\\
d\alpha_3\wedge \beta_3 \wedge \beta_5=0,\\
d\alpha_5\wedge \beta_3 \wedge \beta_5=0,\\
\end{array}
\right.
$$
and $\theta_2(X_3,Y)=\theta_2(X_5,Y)=0$ for any $Y$. Thus there exists linear forms $\gamma_j \in \K\{\alpha_i, i\neq 3,5\}$ with
$$
\left\{
\begin{array}{l}
d\alpha_3=\beta_3\wedge\gamma_1+\beta_5\wedge\gamma_2  \\
d\alpha_5=\beta_3\wedge\gamma_3+\beta_5\wedge\gamma_4,
\end{array}
\right.
$$
and

$$\beta_3\wedge \beta_5\wedge \gamma_2-\beta_3\wedge \beta_5\wedge \gamma_3=0.$$
This implies
$$\gamma_3=\gamma_2 + a \beta_3 + b\beta_5$$
and
$$
\left\{
\begin{array}{l}
d\alpha_3=\beta_3\wedge\gamma_1+\beta_5\wedge\gamma_2  \\
d\alpha_5=\beta_3\wedge\gamma_2+\beta_5\wedge (b\beta_3+\gamma_4),
d\alpha_i=0, \ i \neq 3,5.
\end{array}
\right.
$$
We deduce
\begin{theorem}
Any $2p$-dimensional symplectic $2$-step nilpotent Lie algebras with characteristic sequence $(2,2,1,\cdots,1)$ and $p \geq 4$ is decomposable and isomorphic to   $\n_7^{124} \oplus \K^{2p-7}$,  $\n_6^{19} \oplus \K^{2p-6}$,   $\n_6^{20} \oplus \K^{2p-6}$, $\n_5^{5} \oplus \K^{2p-5}$, $\n_3^{1}\oplus \n_3^1 \oplus \K^{2p-3}$ for $\K=\C$ or $\R$ and $\n_6^{20,1} \oplus \R^{2p-6}$.
\end{theorem}

\subsection{Symplectic Lie algebras in $\mathcal{F}^{2p,2}_{k,2p-2k}$}

Lie algebras of the set $\mathcal{F}^{2p,2}_{k,2p-2k}$ are of characteristic sequence $(2,\cdots,2,1,\cdots,1)$  where $2$ appears $k$-times.  Let $\g$ be in $\mathcal{F}^{2p,2}_{k,2p-2k}$. We assume that $\g$ admits $2p-(k+1)$ generators. The others cases could be studied in a same way. We consider a basis adapted to the characteristic sequence, the dual basis satisfying the Maurer-Cartan equations:
$$
\left\{
\begin{array}{l}
d\alpha_3=\alpha_1\wedge\alpha_2+\sum a_{ij}^3\alpha_i\wedge\alpha_j, \ \ 2 \leq i <j \leq 2p, \ i,j \neq 3, 5,\cdots,2k+1,  \\
d\alpha_5=\alpha_1\wedge\alpha_4+\sum a_{ij}^5\alpha_i\wedge\alpha_j, \ \ 2 \leq i <j \leq 2p, \ i,j \neq 3, 5,\cdots,2k+1,\\
\cdots \\
d\alpha_{2k+1}=\alpha_1\wedge\alpha_{2k}+\sum a_{ij}^{2k}\alpha_i\wedge\alpha_j, \ \ 2 \leq i <j \leq 2p, \ i,j \neq 3, 5,\cdots,2k+1,\\
d\alpha_i=0, \ i \neq 3,5.\cdots,2k+1.
\end{array}
\right.
$$
Let $\theta$ be a symplectic form on $\g$. It can be decomposed in $\theta=\theta_1+\theta_2$ with $\theta_1=\alpha_3 \wedge \beta_3+\alpha_5 \wedge \beta_5+\cdots +\alpha_{2k+1}\wedge \beta_{2k+1}$, $d\theta_1=0$, $\theta_2(X_3,Y)=\theta_2(X_5,Y)=\cdots =\theta_2(X_{2k+1},Y)=0$ for any $Y$ and $\theta_2^k\neq 0$.  This implies, since by hypothesis $d\beta_i=0$,
$$
\left\{
\begin{array}{l}
d\alpha_3\wedge \beta_3+d\alpha_5 \wedge \beta_5+\cdots+d\alpha_{2k+1} \wedge \beta_{2k+1}=0,\\
d\alpha_3\wedge \beta_3 \wedge \beta_5 \wedge\cdots\wedge \beta_{2k+1} =0,\\
d\alpha_5\wedge \beta_3 \wedge \beta_5\wedge\cdots\wedge \beta_{2k+1}=0,\\
\cdots \\
d\alpha_{2k+1}\wedge \beta_3 \wedge \beta_5\wedge\cdots\wedge \beta_{2k+1}=0,
\end{array}
\right.
$$
 Thus there exists linear forms $\gamma_{i,j} \in \K\{\alpha_i, i\neq 3,5,\cdots,2k+1\}$ with
$$
\left\{
\begin{array}{l}
d\alpha_3=\beta_3\wedge\gamma_{3,1}+\beta_5\wedge\gamma_{3,2}+\cdots+ \beta_{2k+1}\wedge\gamma_{3,k}, \\
d\alpha_5=\beta_3\wedge\gamma_{5,1}+\beta_5\wedge\gamma_{5,2}+\cdots+ \beta_{2k+1}\wedge\gamma_{5,k}, \\
\cdots\\
d\alpha_{2k+1}=\beta_3\wedge\gamma_{2k+1,1}+\beta_5\wedge\gamma_{2k+1,2}+\cdots+ \beta_{2k+1}\wedge\gamma_{2k+1,k}, \\
\end{array}
\right.
$$
and
$$\gamma_{2i+1,j}=\gamma_{2j+1,i} \ {\text{\rm mod}} \ \beta_3,\cdots,\beta_{2k+1}.$$
We deduce
\begin{theorem}
Let $\g$ be a $2$-step nilpotent symplectic Lie algebra with characteristic sequence $(2,\cdots,2,1,\cdots,1)$ and dimension $2p=2k+l$ where $k$ is the number of $2$ in the characteristic sequence. If $\ds \dim \g=2p > \frac{k(k+5)}{2}$, then $\g$ is decomposable. It is a direct product of a $\ds  \frac{k(k+5)}{2}$-dimensional Lie algebra by an abelian ideal.
\end{theorem}

\noindent{\bf Remark.} Let $\g$ be a $2$-step nilpotent Lie algebra with  $c(\g)=(2,\cdots,2,1,\cdots,1)$ and dimension $2p=2k+l$ where $k$ is the number of $2$ in the characteristic sequence. We assume that  $g=\h \oplus \mathcal{A}_s$ where $\mathcal{A}_s$ is an abelian ideal. Then $c(\h)=(2,\cdots,2,1,\cdots,1)$ and 
$k$ is always the number of $2$ in this sequence. As consequence, if $\g \in \mathcal{F}^{2p,2}_{k,2p-2k}$ is symplectic and indecomposable, then $\dim \g \leq \ds \frac{k(k+5)}{2}$.

\subsection{Symplectic Lie algebras in $\mathcal{F}^{2p,2}_{p-1,2}$}
From \cite{GRKegel}, the set $\mathcal{F}^{2p,2}_{p-1,2}$ is the union of two algebraic components, each one being characterized by the number of generators. We consider here the component $\mathcal{C}_1(\mathcal{F}^{2p,2}_{p-1,2})$ because we have seen that in dimension $8$ some Lie algebras of this family are not symplectic. A Lie algebra of $\mathcal{C}_1(\mathcal{F}^{2p,2}_{p-1,2})$ is isomorphic to 
$$
\left\{
\begin{array}{l}
[X_1,X_{2i}]=X_{2i+1}, \ i=1,\cdots p-1,  \\
\lbrack X_{2i},X_{2j}\rbrack =\displaystyle \sum_{k=1}^{p-1}C_{2i,2j}^{2k+1}X_{2k+1}, \ 1\leq i<j \leq p.
\end{array}
\right.
$$
If $\g$ is symplectic, any symplectic form $\theta $ on $\g$ satisfies
$$\theta(X_i,X_j)=0$$
for $i,j=3,5,\cdots,2p-1.$  Consider the decomposition  $\theta=\theta_1+\theta_2$ with 
$$\theta_1=\alpha_3 \wedge \beta_3 + \alpha_5 \wedge \beta_5 +\cdots+ \alpha_{2p-1} \wedge \beta_{2p-1}$$
and $\beta_i(X_j)=0$ for $i,j=3,5,\cdots,2p-1$ and $\theta_2=\sum a_{i,j} \alpha_i \wedge \alpha_j$ with $i,j \in I=\{1,2,4,\cdots,2p\}$. Since $d\alpha_i=0$ for any $i \in I$, we have $d\theta_2=0$. Then $\theta$ symplectic implies $d\theta_1=0$.
This is equivalent to
$$d\alpha_3 \wedge \beta_3 +d\alpha_5 \wedge \beta_5 +\cdots+d\alpha_{2p-1} \wedge \beta_{2p-1}=0.$$
Similarly to the case of dimension $8$, 
 we  characterize the Lie algebras such that any closed $2$-forms is not of maximal rank, that is satisfies $\theta^p=0$. This last condition is equivalent to $\theta_1^{p-1} =0,$ that is, $ \beta_3\wedge \beta_5  \wedge\cdots \wedge \beta_{2p-1} =0$.  The linear forms $\beta_i, i=3,5,\cdots,2p-1$ are not independent 
and there exists  a change of basis respecting the characteristic sequence giving 
$\beta_{2p-1}=0$. Then $\g$ is not symplectic if the condition $d\theta_1 =0$ implies $\beta_{2p-1}=0$.  A  computation similar to the one in section 4.2 shows that $\alpha_{2p-1}$ is of maximal Cartan Class and $\g$ is isomorphic to
$$
\left\{
\begin{array}{l}
d\alpha_3=\alpha_1 \wedge \alpha_2, \\
d\alpha_5=\alpha_1 \wedge \alpha_4, \\
\cdots\\
d\alpha_{2p-3}= \alpha_1 \wedge \alpha_{2p-4},\\
d\alpha_{2p-1}=\alpha_1 \wedge \alpha_{2p-2}+\sum_{1\leq i < j \leq p} C^{2k+1}_{2i,2j}\alpha_{2i}\wedge\alpha_{2j}, \\
d\alpha_{2p}=0.
\end{array}
\right.
$$
Such a Lie algebra doesn't admit a general symplectic form (that is without the condition $\beta_i(X_1)=0$) if and only if $C^{2p-1}_{2p-2,2p} \neq 0.$
In this case, since $\alpha_{2p-1}$ is of maximal Cartan Class and $C^{2p-1}_{2p-2,2p} \neq 0$, we can reduce the form
$d\alpha_{2p-1}$ to
$$d\alpha_{2p-1}=\alpha_1 \wedge \alpha_{2p-2}+\alpha_{2p-2}\wedge \alpha_{2p}+ \alpha_{2}\wedge \alpha_{4}+\cdots+\alpha_{2p-6}\wedge \alpha_{2p-4}$$
if $p$ is even or
$$d\alpha_{2p-1}=\alpha_1 \wedge \alpha_{2p-2}+\alpha_{2p-2}\wedge \alpha_{2p}+ \alpha_{2}\wedge \alpha_{4}+\cdots+\alpha_{2p-8}\wedge \alpha_{2p-6}$$
if $p$ is odd.

\end{document}